\documentclass[11pt]{article}
\usepackage{amsmath,amssymb,amsthm, mathrsfs}

\begin{document}

\renewcommand{\baselinestretch}{1.3}
\renewcommand{\arraystretch}{1.3}

\begin{center}{\bf\LARGE  Complete  Rank Theorem in Advanced Calculus and  Frobenius Theorem in Banach Space}\\
\vskip 0.5cm Ma Jipu$^{1,2}$
\end{center}

{\bf Abstract}\quad  The application of generalized inverses is usually neglected in pure mathematical research.
However, it is very effective for
 this paper. The famous matrix rank theorem  was established by R. Penrose in 1955.
   We expand the theorem from the case of matrices to that of  operators.
 Then a modern perturbation analysis of generalized inverses is built. Hereby, we find and prove a complete
  rank theorem in advanced calculus. So  a  complete answer to the rank theorem problem presented by M. S. Berger is
  given. Let $\Lambda$ be an open set in a Banach space $E$ and $x_0$ be a point in $\Lambda$.
  We consider the family of subspaces
  $\mathcal{F}$$=\{M(x)\}_{x\in\Lambda}$,  especially, where  the
  dimension of $M(x)$ may be infinity,   and investigate the necessary and sufficient
condition for $\cal{F}$ being $c^1$  integrable at  $x_0$. For this
we introduce the concept of the  co-final set $J(x_0,E_*)$ of
$\cal{F}$ at $x_0$. Then applying  the co-final set and the
perturbation analysis of generalized inverses we prove  the
Frobenius theorem in Banach space in the proof of which the used
vector field and flow theory  are avoided.  The co-final set is
essential to the Frobenius theorem. When $J(x_0,E_*)$ is trivial,
the theorem reduces to an initial value problem for a differential equation in Banach space.
Let $B(E,F)$ be the set of all bounded linear
operators from  $E$ into another Banach space $F$,
$\Lambda=B(E,F)\setminus\{0\}$ and $M(X)=\{T\in B(E,F):TN(X)\subset
R(X)\}$ for any $X\in\Lambda$. In 1989, V. Cafagna introduced a
geometrical method for some partial differential equations and  the
family of subspaces $\mathcal{F}=\{M(X)\}_{x\in\Lambda}$. Let
$\Phi_{1,1}$ be the set of all Fredholm operators $T$ with
dim$N(T)=$codim$R(T)=1$.  It is essential  to this method that
$\Phi_{1,1}$ is a smooth submanifold in $B(E,F)$ and tangent to
$M(X)$ at any $X\in\Phi_{1,1}$.  The co-final set of the  family $\cal{F}$ at
$A\in\Lambda$ is nontrivial in general. Applying the perturbation
analysis of generalized inverses, the above property for
$\Phi_{1,1}$ is expanded to  wider classes of operators
$\Phi_{m,n}$, semi-Fredholm operators $\Phi_{m,\infty}$ and
$\Phi_{\infty,n}$. It seems to be useful for further developing the
method introduced by V. Cafagna. \vskip1cm



\newpage
\vskip 0.2cm\begin{center}{\bf0\quad Introduction} \end{center}
\vskip 0.2cm

Let $E,F$ be two Banach spaces, $U$ an open set in $E$ and $f$ a
nonlinear $c^1$ map from $U$ into $F$. Recall that $f$ is said to be
locally conjugate to $f'(x_0)$ (Fr$\acute{\rm e}$chet derivative of
$f$ at $x_0)$ near $x_0\in U$ provided there exist a neighborhood
$U_0(\subset U)$ at $x_0$, a neighborhood $V_0(\subset F)$ at 0, and
diffeomorphisms $\varphi:U_0\rightarrow \varphi(U_0)$ with
 $\varphi(x_0)=0$,  $\varphi'(x_0)=I_E$ and $\psi:V_0\rightarrow
 \psi(V_0)$ with $\psi(0)=y_0(=f(x_0)),$  $\psi'(0)=I_F$, such
 that
 $$f(x)=(\psi\circ f'(x_0)\circ\varphi)(x),\quad\quad\forall x\in
 U_0,\eqno(*)$$
 where $I_E$ and $I_F$ denote  the identities on $E$ and $F$, respectively.

 It is well known that if one of the following conditions holds:
 $N(f'(x_0))=\{0\}$
 and $R(f'(x_0))$ splits in $F; R(f'(x_0))=F$,
 $N(f'(x_0))$ splits in $E$; and rank$f'(x)=$
 rank$f'(x_0)<\infty$
 near $x_0$, then $f$ is locally conjugate to $f'(x_0)$ near $x_0$.
 These are important basic theorems in non-linear functional
 analysis (for details, see [Zie] and [Abr]).
 By saying the  rank theorem problem one means what property of $f$ (more
 general than the three conditions above) ensures the equality $(*)$
 holds.
 In [Beg], M. S. Berger points out that it is not yet known whether the rank
 theorem in advanced calculus holds even if $f$ is a Fredholm map.
 It is not difficult to observe that each of the three conditions
 above ensures that $f'(x)$ near $x_0$ has a  generalized inverse
 $f'(x)^+$ (See Section 1) satisfying
 $$\lim\limits_{x\rightarrow x_0}f'(x)^+=f'(x_0)^+.\eqno(**)$$
 So we have in  mind to seek the condition for which  the equality
 (**) holds,  and try to find the answer to the rank theorem
 problem.

Let $X$ be a topological space, and $B(E, F)$ be the set of all bounded linear operators from  $E$ into  $F$.
Suppose that the operator valued map $T_x: X\mapsto B(E, F)$ is continuous at $x_0\in X$, and $T_{x_0}$ is
 double splitting.  In Section 1, we will introduce the concept of a locally fine point of $T_x$ and
prove the following operator rank theorem: for any generalized
inverse $T_0^+$ of $T_{x_0}$, there exists a neighborhood $U_0$ at
$x_0$ such that $T_x$ for any $x\in U_0$ has a generalized inverse
$T_x^+$ and $T_x^+\rightarrow T_0^+$ as $x\rightarrow x_0$ if and only if  $x_0$ is a
locally fine point of $T_x$. Now it
is obvious that  the equality $(**)$ holds if and only if $x_0$ is a
locally fine point of $f'(x)$. As expected, the following complete
rank theorem in advanced calculus will be proved in this section:
$f$ is  locally conjugate to $f'(x_0)$, i.e., the equality $(*)$
holds, if and only if  $x_0$ is a locally fine point of $f'(x)$. By
Theorem 1.9, one observes that the complete rank theorem in advanced
calculus expands widely the three results known well. Therefore  the
rank
 theorem problem presented by M. S. Berger
 has  a complete answer.  Theorems 1.10 and 1.11 are essential to the form of Frobenius theorem
 in Banach space and smooth differential structure for some smooth submanifolds, respectively.
 In Section 2,
 we consider the family of
 subspaces in $E$ as, $\mathcal{F}$ $=\{M(x):x\in\Lambda\}$, where
 $\Lambda$ is an open set in $E$, and $M(x)$ is a subspace in $E$,
 especially, where  the
  dimension of $M(x)$ may be infinity.  Recall that $\mathcal{F}$ is said to
 be $c^1$ integrable at $x_0$ provided there exist a neighborhood
 $U_0$ at $x_0$ and a $c^1$-submanifold $S$ in $U_0$, such that
 $x_0\in S$ and
 $$M(x)=\{\dot{c}(0):\quad{\rm  for\ all}\quad c^1{\rm -curve}\ c(t)\subset S\ {\rm
 with}\ c(0)=x\},$$
 $\forall x\in S$, i.e., $S$ is tangent to $M(x)$ at $x\in
 S$. We introduce the co-final set $J(x_0,E_*)$ of
 $\mathcal{F}$ at $x_0\in\Lambda$  and prove if $\mathcal{F}$
 is $c^1$ integrable at $x_0$, say that $S$ is the $c^1$ integral
 submanifold in $E$, then there exists a neighborhood $U_0$ at
 $x_0$, such that $J(x_0,E_*)\supset S\cap U_0$. Moreover, we
 have the coordinate expression of $M(x)$ as,
 $M(x)=\{e+\alpha e:\forall e\in M_0\}$ for all $x\in J(x_0,E_*)$,
 where $M_0=M(x_0)$ and $\alpha\in B(M_0,E_*).$ Therefore we obtain the
 necessary and sufficient condition for $\mathcal{F}$ being $c^1$
 integrable at $x_0$, i.e., the  Frobenius theorem in  Banach space.
 Thanks to the co-final set $J(x_0,E_*)$, the used vector field and
 flow theory are avoided in the proof of the theorem. The co-final
 set is essential to the theorem. When $J(x_0,E_*)$ is trivial, the
 theorem reduces to an initial value problem for a differential equation in Banach space.
  Let $\Lambda=B(E,F)\setminus\{0\}$ and $M(X)=\{T\in
 B(E,F):TN(X)\subset R(X)\}$ for any $X\in\Lambda$. In section 3, we
 investigate the family of subspaces
 $\mathcal{F}=\{M(X)\}_{X\in\Lambda}$. In fact, V. Cafagna  introduces the
 geometrical method for some partial differential equations in
 [Caf],  and presents the following property:
 let
  $\Phi_{1,1}$ be the set of all
 Fredholm operators $T$ with dim$N(T)=$co-dim$R(T)=1$, then  $\Phi_{1,1}$ is a  smooth
 submanifold in $B(E,F)$ and tangent to $M(X)$ at any
 $X\in\Phi_{1,1}$, which is the key point   to the method. The
 co-final set of $\cal{F}$ at the double splitting $A\in\Lambda$, in
 general, is non-trivial. In this case,in order to get the integral submanifold
 $S$ of $\cal{F}$ at $A$, we try its co-final set. Let $\Phi_{m,n},\Phi_{m,\infty}(\Phi_{\infty,n})$ and
 $F_k$ denote the sets of all Fredholm operators $T$ with dim$N(T)=m$
 and co-dim$R(T)=n$,  of semi-Fredhlom  operators $T$ with dim$N(T)=m$
 and co-dim$R(T)=\infty$(dim$N(T)=\infty$ and co-dim$R(T)=n)$, and
 of
 operators $T$ with rank$T=k<\infty$, respectively. Applying the modern
 perturbation analysis of generalized inverses we prove the
 following results: suppose that $\Phi$ is any one of
 $\Phi_{m,n},\Phi_{m,\infty} (\Phi_{\infty,n})$ and
 $F_k$ , then $\Phi$ is a smooth submanifold in   $B(F, E)$,
 and   tangent to $M(X)$ at any
 $X\in \Phi$. Obviously, these results expand the property for $\Phi_{1,1}$  proposed
by V.  Cafagna to  wider classes of operators. It  seems to be
 useful for further developing  this method (see [An]).

 \vskip 0.2cm
\begin{center}{\bf 1\quad   Complete Rank Theorem and Modern Perturbation Analysis of Generalized Inverses} \end{center}\vskip 0.2cm

Recall that $A^+\in B(F, E)$ is said to be a generalized inverse of
$A\in B(E,F)$\ provided $A^+AA^+=A^+$ and $A=AA^+A$, (when both $E$
and $F$ are Hilbert spaces, $A^+$ is said to be M. -P. inverse of
$A$ provided $(AA^+)^*=AA^+$ and $(A^+A)^*=A^+A$, and it  is
unique); $A\in B(E,F)$ is said to be double splitting if $R(A)$ is
closed, and there exist closed subspaces $R^+$ in $E$ and $N^+$ in
$F$ such that $E=N(A)\oplus R^+$ and $F=R(A)\oplus N^+$,
respectively. It is well known that $A$ has a generalized inverse
$A^+\in B(F,E)$ if and only if $A$ is double splitting.

The classical perturbation analysis of the generalized inverse is
marked with the matrix rank theorem by R. Penrose. Let $A$ and
$\Delta A$ belong to $B(\mathbf{R}^n),$ and both  $A^+$ and
$(A+\Delta A)^+$ be M. -P. inverses of $A$ and $(A+\Delta A),$
respectively. The theorem says that $\lim\limits_{\Delta
A\rightarrow o}(A+\Delta A)^+=A^+$ if and only if rank$(A+\Delta
A)=$ rank $A$ for $\Delta A$ small enough. The condition in the
theorem ensures that the computing of $A^+$ is stable. Hereafter one
 will try to expand the theorem to the case of operators in $B(E,F)$. This
is very difficult because of no fixed one of generalized inverses
like M. -P. inverse for every double splitting operator in
$B(E,F)$. Thanks to a study  of rank theorem problem proposed  by
M. S. Berger in [Beg], we  find the   main concept and theorems in
modern perturbation analysis of generalized inverses, and obtain a
complete rank theorem. First we will work on perturbation analysis  of
generalized inverse. Let $A\in B(E,F)$ be double splitting,
$A\not=0,$ and $A^+$ be a generalized inverse of $A$. Write
$V(A,A^+)=\{T\in
B(E,F):\|T-A\|<\|A^+\|^{-1}\},C_A(A^+,T)=I_F+(T-A)A^+,$ and
$D_A(A^+,T)=I_E+A^+(T-A).$ Then we have

{\bf Theorem 1.1}\quad {\it The following conditions for $T\in
V(A,A^+)$ are equivalent} :

(i) $R(T)\cap N(A^+)=\{0\}$;

(ii) $B=A^+C_A(A^+,T)=D_A(A^+,T)A^+$ {\it is the  generalized
inverse of $T$ with $R(B)=R(A^+)$ and $N(B)=N(A^+)$};

(iii) $R(T)\oplus N(A^+)=F;$

(iv) $N(T)\oplus R(A^+)=E$;

(v) $(I_E-A^+A)N(T)=N(A)$;

(vi) $C^{-1}_A(A^+,T)TN(A)\subset R(A)$;

(vii) $R(C^{-1}_A(A^+,T)T)\subset R(A)$.

{\bf Proof}\quad  We are going to show ${\rm(vi)\Leftrightarrow(ii)}$. Note that M. Z. Nashed and  X. C. Chen have proved
(vi)$\Rightarrow$ (ii) in [N-C]. Indeed, the inverse relation holds
too. Obviously, $C^{-1}_A(A^+,T)TA^+A=A$ for all $T\in V(A,A^+)$.
 Hereby it follows that $$BTB-B=0\quad{\rm and}\quad
TBT-T=-(I_F-AA^+)C^{-1}(A^+,T)T\eqno (1.1)$$ for any $T\in
V(A,A^+)$. So(vi)$\Leftrightarrow$(ii).

Evidently,
$$\begin{array}{rllr}
C^{-1}_A(A^+,T)Th&=&C^{-1}_A(A^+,T)TA^+Ah+C^{-1}_A(A^+,T)T(I_E-A^+A)h\\
&=&Ah+C^{-1}_A(A^+,T)T(I_E-A^+A)h,\quad\forall h\in E,\end{array}$$
so that (vi)$\Leftrightarrow$(vii).

Go to the  claim (v)$\Leftrightarrow$(vi). Assume that (vi) holds.
We have for any $h\in N(A)$, there exists $g\in R(A^+)$ such that
$C^{-1}_A(A^+,T)Th=Ag$. Note $A^+Ag=g,$ then
$C^{-1}_A(A^+,T)Th=Ag=C^{-1}_A(A^+,T)TA^+Ag=C^{-1}_A(A^+,T)Tg.$ So
$h-g\in N(T)$ and satisfies $(I_E-A^+A)(h-g)=h$. This shows that (v)
holds. Conversely, assume that (v) holds. Then we have for any $h\in
N(A)$, there exists $g\in N(T)$ such that $h=(I_E-A^+A)g,$ and so,
$C^{-1}_A(A^+,T)Th=C^{-1}_A(A^+,T)TA^+Ag=Ag\in R(A)$, i.e., (vi) holds. This shows
(v)$\Leftrightarrow$ (vi).

Go to the claim (i)$\Leftrightarrow$(ii). Obviously,
(ii)$\Rightarrow$(i). Conversely, assume that (i) holds.
 By (1.1), $BTB-B=0$
and $R(TBT-T)\subset N(A^+).$  Obviously $R(TBT-T)\subset R(T)$.  By
the assumption $R(T)\cap N(A^+)=\{0\}$, one concludes that $B$ is
the generalized inverse of $T$, this shows (i)$\Rightarrow$(ii). So
(i)$\Leftrightarrow$(ii).

Go to show (i)$\Leftrightarrow$(iii). Obviously,
(iii)$\Rightarrow$(i). Conversely, assume that (i) holds. Then $B$
is the generalized inverse of $T$ with $N(B)=N(A^+)$ because of the
equivalence of (i) and (ii). So (i)$\Leftrightarrow$(iii).

Finally
 go to show (i)$\Leftrightarrow$(iv). Assume that (i) holds.  Since
 (i)$\Rightarrow$(ii), one can conclude
(iv) holds.
 Conversely, assume that (iv) holds. Then
 $N(A)=(I_E-A^+A)E=(I_E-A^+A)N(T)$ i.e., (v) holds. Hence (i) holds
 since  (v)
 $\Leftrightarrow$(vi)$\Leftrightarrow$(ii)$\Leftrightarrow$(i).
 This shows (i)$\Leftrightarrow$(iv).

 So far we have proved the following relations
 (i)$\Leftrightarrow$(ii),(i)$\Leftrightarrow$(iii),(i)$\Leftrightarrow$(iv),
 (ii)$\Leftrightarrow$(vi),
 (v)$\Leftrightarrow$(vi), and (vi)$\Leftrightarrow$(vii). Thus the
 proof of the
 theorem is completed.\quad$\Box$

  (For more information, see [H-M]).

 The condition (i) looks strange; however, it may imply  many
 known conditions  in analysis as showing in the next theorem.

 Let $F_k=\{T\in B(E,F):$ rank$ T=k<\infty\}$, and  $\Phi_{m,n}=\{T\in B(E,F)$: dim $N(T)=m<\infty$ and
  codim$R(T)=n<\infty\}$.
 Let  $\Phi_{m,\infty}$ be the set of all semi-Fredholm operators
 $T$ with dim$N(T)=m$ and codim$R(T)=\infty$, and $\Phi_{\infty,n}$
 be the set of all semi-Fredholm operators $T$ with dim$N(T)=\infty$
 and codim$R(T)=n$. We have

 {\bf Theorem 1.2}\quad {\it Assume that  $A$  belongs to  any one of $F_k,
 \Phi_{m,n},\Phi_{m,\infty}$ and $\Phi_{\infty,n}.$ Then  the condition
 $R(T)\cap N(A^+)=\{0\}$ for $T\in V(A,A^+)$ holds if and only if
 $T$ belongs to the  corresponding class for  $A$.}

 {\bf Proof}\quad Assume that the condition $R(T)\cap
 N(A^+)=\{0\}$ holds for $T\in V(A,A^+)$. By (ii) in Theorem 1.1,
 $T$ has a generalized inverse $B$ with $N(B)=N(A^+)$ and
$R(B)=R(A^+)$,
 so that
 $$N(T)\oplus R(A^+)=E=N(A)\oplus R(A^+)$$
 and
$$R(T)\oplus N(A^+)=F=R(A)\oplus N(A^+).$$
 Hereby one observes   that when $A$
 belongs to any one of  $F_k,\Phi_{m,n},\Phi_{m,\infty}$ and $\Phi_{\infty,n},T$
  belongs to the   corresponding  class for  $A$.

 Conversely, assume that $T$  belong to  any one of
 $F_k,\Phi_{m,n},\Phi_{m,\infty}$ and $\Phi_{\infty,n}$.
By (1.1),
 $B=A^+C^{-1}_A(A^+,T)=D^{-1}_A(A^+,T)A^+$ satisfies
$BTB=B.$  Thus $B$ and $T$ bear two projections $P_1=BT$ and
$P_2=TB$.
 Indeed, $P^2_1=BTBT=BT=P_1,$ and $P^2_2=TBTB=TB=P_2$.
Note $N(B)=N(A^+)$. Clearly, $N(P_1)=N(T)+\{e\in R(T^+): Te\in
N(A^+)\}$,  where by the assumption, $T$ has a generalized inverse
$T^+$. We next claim
 $$R(P_1)=R(A^+), R(P_2)=R(TA^+),\quad{\rm and}\quad
 N(P_2)=N(A^+).$$

Obviously, $R(P_1)\subset R(B)= R(A^+)$. On the  other hand, by
$C^{-1}_A(A^+,T)TA^+=AA^+$ we get
 $$P_1A^+e=A^+C^{-1}_A(A^+,T)TA^+e=A^+AA^+e=A^+e,\quad\forall e\in
 E.$$
 So
 $R(P_1)=R(A^+)$.

 Obviously,
 $$R(P_2)=R(TA^+C^{-1}_A(A^+,T))=R(TA^+).$$
 $$N(P_2)=N(TD^{-1}_A(A^+,T)A^+)\supset N(A^+),$$
and
 $$Be=BTBe=BP_2e=0,\quad\forall c\in N(P_2).$$
 So, $N(P_2)=N(A^+)$ because of $B=D^{-1}_A(A^+,T)A^{-1}.$

 Thus we have
 $$F=R(P_2)\oplus N(P_2)=R(TA^+)\oplus N(A^+)\eqno(1.2)$$
 and
 $$E=R(P_1)\oplus N(P_1)=R(A^+)\oplus E_*\oplus N(T),$$
 where
 $E_*=\{e\in R(T^+):Te\in N(A^+)\}.$

 We now are in the position to end the proof.

 Assume that $T\in V(A,A^+)$ satisfies rank$T=$rank$A<\infty$.
 Obviously $T$ is double splitting. Let $T^+$ be an arbitrary
 generalized inverse of $T$, then $E=N(T)\oplus R(T^+)$. Note
 dim$R(T^+)=$rank$T=$rank$A={\rm dim}R(A^+)<\infty.$ According to
 $(1.2)$  ${\rm dim}E_*=0$, i.e., $R(T)\cap
 N(A^+)=\{0\}.$

 Assume that $T\in V(A,A^+)$ satisfies codim$R(T)={\rm
 codim}R(A)=n<\infty$, and $N(T)$  is splitting in $E$. Obviously, $T$
 is double splitting. Let $T^+$ be  an arbitrary  generalized
 inverse of $T$, then $F=R(T)\oplus N(T^+)$. According to (1.2)
 $$R(TA^+)\oplus TE_*\oplus N(T^+)=F=
R(TA^+)\oplus N(A^+).$$
 Since dim$N(A^+)={\rm codim}R(A)={\rm
 codim}R(T)={\rm dim}N(T^+)=n<\infty$, one observes dim$TE_*=0$.
 Note that  $E_*\subset R(T^+)$ and $T|_{R(T^+)}$ is invertible in $B(R(T^+), R(T))$,
  then we conclude dim$E_*=0$, i.e.,
 $R(T)\cap N(A^+)=\{0\}.$

 Assume that $T\in V(A,A^+)$ satisfies dim$N(T)={\rm
 dim}N(A)=m<\infty$, and $R(T)$ is splitting in $F$. Obviously, $T$ is
 double splitting. According to (1.2), we have
 $$R(A^+)\oplus E_*\oplus N(T)=R(A^+)\oplus N(A).$$
 Hereby it is easy to observe dim$E_*=0$, i.e.,  $R(T)\cap
 N(A^+)=\{0\}$. So far one can conclude that both  $T$ in
 $V(A,A^+)$ and $A$ belong to the same one of  $F_k,\Phi_{m,n},\Phi_{\infty,n}$
  and $\Phi_{m,\infty}$, then
 $R(T)\cap N(A^+)=\{0\}.$ \quad$\Box
$

 Let $A^\oplus$ be another generalized inverse of $A$,
 $B=A^+AA^\oplus$, and $\delta=\min\{\|A^+\|^{-1},\|B\|^{-1}\}.$
Moreover, we have

 {\bf Theorem 1.3}\quad{\it  Let $S=\{T\in V(A,A^+): R(T)\cap N(A^+)=\{0\}),$
 and
$V_\delta=\{T\in B(E,F):\|T-A\|<\delta\}.$
 Then}
 $$R(T)\cap N(A^\oplus)=\{0\},\quad\forall T\in S\cap V_\delta.$$

 {\bf Proof}\quad It is easy to see that $B$ is the
 generalized inverse of $A$ with $R(B)=R(A^+)$ and
 $N(B)=N(A^\oplus).$ By  the conditions (i) and (iv) in
 Theorem 1.1, $E=N(T)\oplus R(A^+)$, i.e., $E=N(T)\oplus R(B)$ for
 all $T\in S$. Instead of $A^+$ in Theorem 1.1 by $B$, one observes
 that $E=N(T)\oplus R(B)$ is equivalent to $R(T)\cap N(B)=\{0\}$ for
 $T\in V(A,B)$. So $R(T)\cap N(A^\oplus)=\{0\}$ for $T\in S\cap
 V_\delta.$\quad$\Box$

 We now are in the position to discuss modern  perturbation
 analysis of generalized inverses. Consider the  operator valued map
 $T_x$ from a topological space $X$ into $B(E,F)$.

 {\bf Definition 1.1}\quad {\it
 Suppose that $T_x$ is continuous at
 $x_0\in X$, and  that $T_0\equiv T_{x_0}$ is double splitting.
 $x_0$ is said to be a locally fine point of $T_x$
 provided there exist a generalized inverse $T^+_0$ of
 $T_0$ and a neighborhood $U_0$ (dependent on  $T^+_0)$ at $x_0$,
 such that
 $$R(T_x)\cap N(T^+_0)=\{0\},\quad\quad\forall x\in U_0.\eqno(1.3)$$
 For  the locally fine point, we have the  following important theorem.}

 {\bf Theorem 1.4}\quad{\it  The definition of the locally fine point
 $x_0$ of $T_x$ is independent of the choice of the generalized inverse
 $T^+_0$ of $T_0$.}

 {\bf Proof}\quad Assume that (1.3) holds for a generalized inverse
 $T^+_0$ of $T_0$  and a neighborhood $U_0$ at $x_0$. Let
 $T^\oplus_0$ be another generalized inverse of
 $T_0,\delta=\min\{\|T^+_0\|^{-1},\|T^+_0T_0T^\oplus_0\|^{-1}\}$,
 and $V_\delta=\{T\in B(E,F):\|T-T_0\|<\delta\}$. Set $V_0=\{x\in
 U_0:T_x\in V_\delta\}$. Then by Theorem 1.3, $R(T_x)\cap N(T^\oplus
 _0)=\{0\}$ for all $x\in V_0.$
\quad$\Box$

The following theorem expands the famous matrix rank theorem by R.
Penrose to operators in $B(E,F)$.

{\bf Theorem 1.5}\quad (Operator rank theorem)\quad {\it Suppose
that the operator valued map $T_x:X\rightarrow B(E,F)$ is continuous
at $x_0\in X$ and $T_0$ is double splitting. Then the following
conclusion holds for arbitrary  generalized inverse $T^{+}_0$ of
$T_0$: there exists a neighborhood $U_0$ at $x_0$ such that $T_x$
has a generalized inverse $T^+_x$ for $x\in U_0$, and
$\lim\limits_{x\rightarrow x_0}T^+_x=T_0$, if and only if $x_0$ is a
locally fine point of $T_x$.}

{\bf Proof}\quad Assume that $x_0$ is a locally fine point of $T_x$.
Let $U_0=\{x\in X:\|T_x-T^+_0\|<\|T^+_0\|^{-1}\}$ for an arbitrary
generalized inverse $T^+_0$ of $T_0$, and
$T^+_x=T^+_0C^{-1}_{T_0}(T^+_0,T_x)$ for every $x\in U_0$. Since $T_x$
is continuous at  $x_0\in X$, one observes that $U_0$ is a
neighborhood at $x_0$. However,  it is easy to see that $T^+_x$ is the
generalized inverse of $T_x$ for any $x\in U_0$ such that
$\lim\limits_{x\rightarrow x_0}T^+_x=T^+_0$. Conversely, assume that
the following conclusion holds for some generalized inverse $T^+_0$
of $T_0$: there exists  a neighborhood $U_0$ at $x_0$ such  that
$T_x$ has a generalized inverse $T^+_x, \forall x\in U_0,$ and
$\lim\limits_{x\rightarrow x_0}T^+_x=T^+_0$. Consider the
projections, $P_x=I_E-T^+_xT_x$ and $P_0=I_E-T^+_0T_0$. Obviously,
$R(P_x)=N(T_x)$,  $R(P_0)=N(T_0)$, and $P_x\rightarrow P_0$ as
$x\rightarrow x_0$. Let $V_0=\{x\in U_0:\|P_x-P_0\}<1\}\cap\{x\in
U_0:\|T_x-T_0\|<\|T^+_0\|^{-1}\}.$ Then by Problem 4.1 in [Ka],
$P_0N(T_x)=N(T_0)$, i.e., the condition (v) in Theorem 1.1 holds for
all $x\in V_0.$ Thus, by the equivalence of (v) and (i) in Theorem
1.1, we conclude $R(T_x)\cap N(T^+_0)=\{0\}, \forall x\in V_0$. Finally, by Theorem 1.4
we conclude that $x_0$ is a locally fine point of $T_x$.
\quad$\Box$

We now are going to investigate the operator rank theorem in Hilbert
spaces like the matrix rank theorem by R. Penrose. Before working on
this, we introduce the following lemma, which  is  very useful in
the operator theory in Hilbert space, although  it is simple.

{\bf Lemma 1.1}\quad  {\it Suppose that $H$ is a Hilbert space, and
$E_1$ is a closed subspace in $H$. Let $H=E_1\oplus E_2$ while the
relation  $E_1\perp E_2$ is not assumed,  and  $P^{E_2}_{E_1}$ be
the projection corresponding to the decomposition $H=E_1\oplus E_2$.
Then
$$(P^{E_2}_{E_1})^*=P^{E^\perp_1}_{E^\perp_2},$$
where $E^\perp_1$ and $E^\perp_2$ denote the orthogonal complements
of $E_1$ and $E_2$, respectively.}

{\bf Proof}\quad For abbreviation, write $P=P^{E_2}_{E_1}$. It is
clear that $P^*$ is also a projection since $(P^*)^2=(P^2)^*=P^*. $
So in order to show the lemma, it is enough to examine
$R(P^*)=E^\perp_2$ and $N(P^*)=E^\perp_1.$ Obviously
$$H=N(P)\oplus R(P^*)=E_2\oplus R(P^*),$$
and
$$H=R(P)\oplus N(P^*)=E_1\oplus N(P^*),$$
where $\oplus$ denotes the orthogonal direct sum. Then the lemma
follows.\quad$\Box$

With Lemma 1.1 we can establish the operator rank theorem
in Hilbert spaces as follows.

{\bf Theorem 1.6}\quad{\it Let $H_1$,  $H_2$ be Hilbert spaces, and
$T_x$ an operator valued map from a topological space $X$ into
$B(H_1,H_2)$. Suppose that $T_x$ is continuous at $x_0\in X$ and
$R(T_0)$ is closed . Then the following conclusion holds: for M.-P.
inverse $T^+_0$ of $T_0$ there exists a neighborhood $U_0$ at
$x_0$ such that $T_x$ has M.-P. inverse $T^+_x$ (i.e., $R(T_x)$
 is closed) for $x\in U_0$ and $\lim\limits_{x\rightarrow
x_0}T^+_x=T^+_0,$ if and only if $x_0$ is a locally fine point of
$T_x$.}

{\bf Proof}\quad Assume that the conclusion of  the theorem holds
for M. -P. inverse $T^+_0$ of $T_0$. By Theorem 1.5 it is immediate
that $x_0$ is a locally fine point of $T_x$. Conversely, assume that
$x_0$ is a locally fine point of $Tx$. By Definition 1.1 there
exists a neighborhood $U_0$ at $x_0$, such that $R(T_x)\cap
N(T^+_0)=\{0\}$ for $x\in U_0$. Without  loss of generality  one can
assume $\|T_x-T_0\|<\|T^+_0\|^{-1}$ for all $x\in U_0$ since $T_x$ is
continuous at $x_0$. Thus, $B_x=T^+_0C^{-1}_{T_0}(T^+_0,T_x)$ is the
generalized inverse of $T_x$ with $R(B_x)=R(T^+_0)$ and
$N(B_x)=N(T^+_0)$ for all $x\in U_0$  because of the  equivalence of
(i) and (ii) in Theorem 1.1 and so, $R(T_x)$ is closed. Let $T^+_x$
be M. -P. inverse of $T_x$. Our goal is to verify
$\lim\limits_{x\rightarrow x_0}T^+_x=T^+_0$. Before working on this
we claim that $\lim\limits_{x\rightarrow x_0}T^+_xT_x=T^+_0T_0$ and
$\lim\limits_{x\rightarrow x_0}T_xT^+_x=T_0T^+_0$. By Lemma 1.1 it
is clear
$$(I_{H_1}-T^*_xB^*_x)=P^{R(T^*_x)}_{N(T_0)},\quad\quad\forall x\in
U_0.\eqno(1.4)$$ In fact, $B_xT_x=P^{N(T_x)}_{R(T^+_0)}$ and so
$(I_{H_1}-T^*_xB^*_x)=(I_{H_1}-(B_xT_x)^*)=(I_{H_1}-P^{N(T_0)}_{R(T^*_x)})=P^{R(T^*_x)}_{N(T_0)}$,
then the equality (1.4) follows from Lemma1.1.

First go to show $\lim\limits_{x\rightarrow x_0}T^+_xT_x=T^+_0T_0.$  We have
\begin{eqnarray*}
\|T^+_xT_x-T^+_0T_0\|&\leq&\|T^+_xT_x(I_{H_1}-T^+_0T_0)\|+\|(I_{H_1}-T^+_xT_x)T^+_0T_0\|\\
&=&\|((I_{H_1}-T^+_0T_0)T^+_xT_x)^*\|+\|(I_{H_1}-T^+_xT_x)T^+_0T_0\|\\
&=&\|(I_{H_1}-T^+_0T_0)T^+_xT_x\|+\|(I_{H_1}-T^+_xT_x)T^+_0T_0\|
\end{eqnarray*}
(Note that $(T^+_xT_x)^*=T^+_xT_x$ and $(T^+_0T_0)^*=T^+_0T_0$ since
both $T^+_x$ and $T^+_0$ are M.-P. inverse), while, by (1.4) and
$R(T^+_x)=R(T^*_x)$,
\begin{eqnarray*}
\|(I_{H_1}-T^+_0T_0)T^+_xT_x\|&=&\|T^+_xT_x-(I_{H_1}-T^*_xB^*_x)T^+_xT_x-T^+_0T_0T^+_xT_x\|\\
&=&\|(I_{H_1}-T^+_0T_0)T^+_xT_x-(I_{H_1}-T^*_xB^*_x)T^+_xT_x\|\\
&\leq&\|(I_{H_1}-T^+_0T_0)-(I_H-T^*_xB^*_x)\|=\|T^*_xB^*_x-T^+_0T_0\|
\end{eqnarray*}
and
\begin{eqnarray*}
\|(I_{H_1}-T^+_xT_x)T^+_0T_0\|&=&\|T^+_0T_0-(I_{H_1}-T^+_xT_x)T^*_xB^*_x-T^+_xT_xT^+_0T_0\|\\
&=&\|(I_{H_1}-T^+_xT_x)T^+_0T_0-(I_{H_1}-T^+_xT_x)T^*_xB^*_x\|\\
&\leq&\|T^+_0T_0-T^*_xB^*_x\|.
\end{eqnarray*}
Since $T^*_xB^*_x\rightarrow T^*_0T^{+*}_0=(T^+_0T_0)^*=T^+_0T_0$ as
$x\rightarrow x_0$, one can conclude $\lim\limits_{x\rightarrow
x_0}T^+_xT_x=T^+_0T_0.$

Similarly
\begin{eqnarray*}
\|T_xT^+_x-T_0T^+_0\|&=&\|(I_{H_2}-T_0T^+_0)T_xT^+_x-T_0T^+_0(I_{H_2}-T_xT^+_x)\|\\
&\leq&\|(I_{H_2}-T_0T^+_0)T_xT^+_x\|+\|((I_{H_2}-T_xT^+_x)T_0T^+_0)^*\|\\
&=&\|(I_{H_2}-T_0T^+_0)T_xT^+_x\|+\|(I_{H_2}-T_xT^+_x)T_0T^+_0\|,
\end{eqnarray*}
while
\begin{eqnarray*}
\|(I_{H_2}-T_0T^+_0)T_xT^+_x\|&=&\|(I_{H_2}-T_0T^+_0)T_xT^+_x-(I_{H_2}-T_xB_x)T_xT^+_x\|\\
&=&\|(T_xB_x-T_0T^+_0)T_xT^+_x\|\leq\|T_xB_x-T_0T^+_0\|
\end{eqnarray*}
and
\begin{eqnarray*}
\|(I_{H_2}-T_xT^+_x)T_0T^+_0\|&=&\|(I_{H_2}-T_xT^+_x)T_0T^+_0-(I_{H_2}-T_xT^+_x)T_xB_x\|\\
&\leq&\|T_0T^+_0-T_xB_x\|.\end{eqnarray*}.

Then it follows $\lim\limits_{x\rightarrow x_0}T_xT^+_x=T_0T^+_0$
from $T_xB_x\rightarrow T_0T^+_0$ as $x\rightarrow x_0$.

Finally, go  to  end the proof. Let $G_x=T^+_xT_xB_x$ and
$Q_x=B_xT_xT^+_x.$ By direct computing
$$G_xT_xQ_x=T^+_xT_xB_xT_xB_xT_xT^+_x=T^+_xT_xB_xT_xT^+_x=T^+_xT^+_x=T^+_x,$$
meanwhile, by the two results proved above, one observes
$$\lim\limits_{x\rightarrow x_0}G_x=T^+_0=\lim\limits_{x\rightarrow
x_0}Q_x.$$ Therefore
$$\lim\limits_{x\rightarrow x_0}T^+_x=\lim\limits_{x\rightarrow
x_0}G_xT_xQ_x=T^+_0T_0T^+_0=T^+_0.$$ The proof ends.\quad$\Box$

This theorem is presented in the dissertation of Q. L. Huang, a post
graduate student of mine; see also [H-M]. Here the complete proof of the theorem is given.

So far, a modern perturbation analysis of generalized inverses is
built.

We now are in the position to establish the rank theorem in advanced calculus.
Let $f$ be a $c^1$ map,
defined in  an open set $U\subset E$ into $F$, and $f'(x_0), x_0\in
U$, be double splitting. By  saying the rank theorem problem we
mean that the  property of $f$ can ensure  the following
conclusion holds: there exist neighborhoods $U_0\subset U$ at $x_0,
V_0\subset F$ at $0$,  diffeomorphisms $\varphi$ from $U_0$ onto
$\varphi(U_0)$ and $\psi$ from $V_0$ onto $\psi(V_0)$ such that $
\varphi(x_0)=0, \varphi'(x_0)=I_E, \psi(0)=f(x_0),\psi'(0)=I,$
 and
$$
f(x)=(\psi\circ f'(x_0)\circ\varphi)(x),\quad\quad\forall x\in
U_0,$$
 i.e.,  $f$ is locally conjugate to $f'(x_0)$ near $x_0$. (For details
 see [Beg]). Since a locally fine point $x_0$ of $f'(x)$ is equivalent
 to $(I_E-T^+_0T_0) N(f'(x))=N(f'(x_0))$ near $x_0$(where $T^+_0$ is
 a generalized inverse of $f'(x_0)$), we have proved indeed the following rank
 theorem in [Ma.1] in 1999's. However, there are many
 typing mistakes in the  argument  of the theorem, one
 can hardly read. After some modifications we rewrite the theorem  as follows.

 {\bf Theorem 1.7}\quad{\it Let $f$  be a $c^1$ map, defined in  an open set $U\subset E$
 into $F$, and $f'(x_0)$ its Fr$\acute{\rm e}$chet derivative  at
 $x_0\in U$.
 Suppose that $f'(x_0)$ is double splitting. If $x_0$ is a locally
 fine point of $f'(x)$, then there exist neighborhoods $U_0$ at
 $x_0, V_0$ at $0$,  diffeomorphisms
 $\varphi:U_0\rightarrow\varphi(U_0)$ and
 $\psi:V_0\rightarrow\psi(V_0$ such that $\varphi(x_0)=0,
 \varphi'(x_0)=I_E,\psi(0)=f(x_0),\psi'(0)=I_F$ and}
 $$f(x)=(\psi\circ f'(x)0)\circ\varphi)(x),\quad\quad\forall x\in
 U_0.$$

 {\bf Proof}\quad It is the key point  to the proof of the theorem to
 determine the diffeomorphisms  $\varphi$ and $\psi$.
 Let $T^+_0$ be a generalized inverse of $T_0=f'(x_0)$. Set
 $$\varphi(x)=T^+_0(f(x)-f(x_0))+(I_E-T^+_0T_0)(x-x_0)$$
 and
 $$
 \psi(y)=(f\circ\varphi^{-1}\circ T^+_0)(y)+(I_F-T_0T^+_0)y.$$
 Obviously, $\varphi'(x_0)=I_E$ and $\varphi(x_0)=0$.
 So there exists an open disk $D^E_{r_0}(x_0)\subset
 E$ with center $x_0$ and radius $r_0$ such that
 $\varphi:D^E_{r_0}(x_0)\rightarrow\varphi(D^E_{r_0}(x_0))$ is a
 diffeomorphism. Let $T_0^+D^F_m(0)\subset\varphi(D^E_{r_0}(x_0))$ where
 $D^F_m(0)\subset F$ is an open disk with center $0$ and radius $m$.
 Then $\psi$ is defined on $D^F_m(0)$. Obviously, $\psi'(0)=I_F$  and
 $\psi(0)=f(x_0)$. So  there exists an open disk in $F$, without loss of generality,
 still write it as  $D^F_m(0)$, such that  $\psi:D^F_m(0)\rightarrow\psi(D^F_m(0))$ is a
 diffeomorphism. Clearly,  both  $T_0\circ\varphi$
 and $\|f'(x)-T_0\|$ are continuous at $x_0$. Hence, there exists an open disk $D^E_{r_1}(x_0)$, $r_1<r_0$,
 such that
 $$T_0\circ\varphi(x)\in D^F_m(0)\ \ \text{and}\ \ \|f'(x)-T_0\|<\|T^+_0\|^{-1},\ \ x\in D^E_{r_1}(x_0).$$
Hence,  $\psi\circ T_0\circ\varphi$ is defined on $D^E_{r_1}(x_0)$.
Since $x_0$ is a locally fine point of $f'(x)$,  we have
 $$N(f'(x))=\varphi'(x)^{-1}N(T_0),\quad\forall x\in
 D^E_{r_1}(x_0),\eqno(1.5)$$
 which is essential to the proof. Indeed,
 $$\varphi'(x)=T^+_0f'(x)+(I_E-T^+_0T_0),\quad\forall x\in
 D^E_{r_0}(x_0)$$
 and so
 $$\varphi'(x)N(f'(x))=(I_E-T^+_0T_0)N(f'(x)).$$
 Then by the  equivalence of (i) and (v) in Theorem 1.1,
 $$N(T_0)=(I_E-T^+_0T_0)N(f'(x))=\varphi'(x)N(f'(x)),\quad\forall
 x\in D^E_{r_1}(x_0),$$
 i.e., the equality (1.5) holds.

 Next go to show that for any disk
 $D^E_l(0)\subset\varphi(D^E_{r_1}(x_0))$, there exists an  open disk
 $D^E_{r_2}(x_0), 0<r_2<r_1$, such that
 $$T^+_0(f(x)-f(x_0))-(1-t)(I_E-T^+_0T_0)(x-x_0)\in
 D^E_l(0),\eqno(1.6)$$
 for all $x\in D^E_{r_2}(x_0)$ and $t\in [0,1].$
Note that
\begin{eqnarray*}
&&\|T^+_0(f(x)-f(x_0))-(1-t)(I_E-T^+_0T_0)(x-x_0)\|\\
&&\quad\quad\leq\|T^+_0(f(x)-f(x_0))\|+\|(I_E-T^+_0T_0)(x-x_0)\|,\end{eqnarray*}
and both  $T^+_0(f(x)-f(x_0))$ and $(I_E-T^+_0T_0)(x-x_0)$ are
continuous and vanish at $x_0$. The conclusion (1.6) is obvious.
We fix some $l$ and write the corresponding $ D^E_{r_2}(x_0)$ as $ D^E_{r}(x_0)$. Applying (1.5) and (1.6) we will prove
$$\begin{array}{rl}
&(f\circ\varphi^{-1})(T^+_0(f(x)-f(x_0))+(I_E-T^+_0T_0)(x-x_0))\\
&\quad\quad=(f\circ\varphi^{-1})(T_0^+(f(x)-f(x_0))),\quad\quad\forall
x\in D^E_r(x_0).\hspace{4cm}(1.7)\end{array}$$
 In view of  (1.6)  consider
$$
\Phi(t,x)=(f\circ\varphi^{-1})(T^+_0(f(x)-f(x_0))-(1-t)(I_E-T^+_0T_0)(x-x_0)),$$
for all  $t\in[0,1]$ and $x\in D^E_r(x_0).$

For abbreviation, write
$x_*=T^+_0(f(x)-f(x_0))-(1-t)(I_E-T^+_0T_0)(x-x_0).$ It is seen directly
$$\frac{d\Phi}{dt}=f'(\varphi^{-1}(x_*))\cdot(\varphi^{-1})'(x_*)\cdot(T^+_0T_0-I_E)(x-x_0).$$
Then by (1.5),  $\dfrac{d\Phi}{dt}=0$ for all $t\in
[0,1]$ and $x\in D^E_r(x_0)$. Thus $\Phi(1,x)=\Phi(0,x)$, i.e., the
equality (1.7) holds. Finally go to prove the theorem. By direct computing
\begin{eqnarray*}
(\psi\circ T_0\circ\varphi)(x)&=&\psi(T_0T^+_0(f(x)-f(x_0)))\\
&=&(f\circ\varphi^{-1})(T^+_0(f(x)-f(x_0))))+(I_F-T_0T^+_0)(T_0T^+_0(f(x)-f(x_0)))\\
&=&(f\circ\varphi^{-1})(T^+_0(f(x)-f(x_0))),\end{eqnarray*} while by
(1.7).
\begin{eqnarray*}
(\psi\circ
T_0\circ\varphi)(x)&=&(f\circ\varphi^{-1})(T^+_0(f(x)-f(x_0))+(I_E-T^+_0T_0)(x-x_0))\\
&=&(f\circ\varphi^{-1})(\varphi(x))=f(x),\quad\quad\forall x\in
D^E_r(x_0).\end{eqnarray*} The theorem is proved.\quad$\Box$

By Theorem 1.5, we further have

{\bf Theorem 1.8}\quad (Complete Rank Theorem)\quad {\it Suppose
that $f:U\subset E\rightarrow F$ is a $c^1$ map, and $f'(x_0), x_0\in
U$, is double splitting. Then $f$ is locally conjugate to $f'(x_0)$
near $x_0$ if and only if $x_0$ is a locally fine point of $f'(x)$.}

{\bf Proof}\quad Assume that $x_0$ is a locally fine point of
$f'(x)$, then by Theorem 1.7, $f$ is  locally conjugate to $f'(x_0)$
near $x_0$.

Assume that $f(x)=(\psi\circ f'(x_0)\circ\varphi)(x),\forall x\in
U_0.$ Let $T^+_0$ be arbitrary one of the generalized inverses of
$T_0=f'(x_0)$. Without loss of generality, one can assume
$\|f'(x)-T_0\|<\|T^+_0\|^{-1},\forall x\in U_0.$ Set
$B_x=\varphi'(x)^{-1}\cdot
T^+_0\cdot\psi'(T_0\circ\varphi(x))^{-1}.$ Since $\varphi'(x_0)=I_E$
and $\psi'(0)=I_F$ one observes $\lim\limits_{x\rightarrow
x_0}B_x=T^+_0$. By direct computing
$$f'(x)=\psi'(T_0\varphi(x))T_6\varphi'(x),\ {\rm and\ so}\ B_xf'(x)B_x=B_x
\ {\rm and}\ f'(x)B_xf'(x)=f'(x)$$ for all $x\in U_0.$
 Thus by Theorem 1.5 we conclude $x_0$ is a locally fine point of $f'(x).$
\quad$\Box$

It is well known that
if any one of the following  conditions holds:

(1) $f'(x)\in F_k$
near $x_0$,

(2) $f'(x_0)\in\Phi_{0,n}(n\leq\infty)$, 
 
(3)$f'(x_0)\in\Phi_{m,0}(m\leq\infty)$

\noindent then $f$ is locally conjugate
to $f'(x_0)$ near $x_0$ (see [Zei]). By Theorem 1.2 we have

{\bf Theorem 1.9}\quad{\it Let $f'(x_0)$ belong to  any  one  of the
following classes: $F_k, \Phi_{m,n}, \Phi_{\infty,n}\\(n<\infty)$ and
$\Phi_{m,\infty}(m<\infty)$. Then $f$ is locally conjugate  to
$f'(x_0)$ near $x_0$, if and only if $f'(x)$ near $x_0$ belongs to
 the class corresponding to $f'(x_0)$.}

Obviously, the classes of operators indicated in Theorem 1.9 contain
properly the three classes above. So,  Theorem 1.7  and Theorem 1.8
give  the  complete answer to the rank theorem problem proposed by
M. S. Berger in [Beg].

Besides Theorems 1.5, 1.6, 1.8 and 1.9 the concept of a locally fine
point leads to expansions of the concepts of regular point, regular
value, and transversility, so that it  bears the  generalized
preimage  and then the transversility theorem by R. Thom (see [Ma2]
and [Ma3]). From  this point of  view, the concept of a locally fine
point should be a nice  mathematical concept.

Also, in the next sections 2 and 3, we will need the following
theorems in perturbation analysis of generalized inverses.

{\bf Theorem 1.10}\quad {\it Suppose that $E_0$ and $E_1$ are two
closed subspaces in a Banach space $E$ with a common complement
$E_*$. Then there exists a unique operator $\alpha\in B(E_0,E_*)$
such that
$$E_1=\{e+\alpha e:\forall e\in E_0\}.$$
Conversely, $E_1=\{e+\alpha e:\forall e\in E_0\}$ for any $\alpha\in
B(E_0,E_*)$ is a closed subspace satisfying $E_1\oplus E_*=E$.}

{\bf Proof}\quad First go to show that $\alpha$ is unique for which
$E_1=\{e+\alpha e:\forall e\in E_0\}$. If $\{e+\alpha e:\forall e\in
E_0\}=\{e+\alpha_1e:\forall e\in E_0\}$, then for any $e\in E_0$
there exists $e_1\in E_0$ such that $(e-e_1)+(\alpha e-\alpha_1
e_1)=0$, and so, $e=e_1$  and $\alpha e=\alpha_1e, $ i.e.,
$\alpha_1=\alpha$.  This says that $\alpha$ is unique. We now claim
that $\alpha=\left.P^{E_0}_{E_*}P^{E_*}_{E_1}\right|_{E_0}\in
B(E_0,E_*)$
 fulfils $E_1=\{e+\alpha e:\forall e\in E_0\}$.

Obviously,
$$P^{E_*}_{E_0}P^{E_*}_{E_1}e=P^{E_*}_{E_0}(P^{E_*}_{E_1}e+P^{E_1}_{E_*}e)=P^{E_*}_{E_0}e=e,\quad\quad\forall
e\in E_0,$$ and
$$P^{E_*}_{E_1}P^{E_*}_{E_0}e=P^{E_*}_{E_1}(P^{E_*}_{E_0}e+P^{E_0}_{E_*}e)=P^{E_*}_{E_1}e=e,\quad\quad\forall
e\in E_1.$$ Hereby let $\alpha=P^{E_0}_{E_*}P^{E_*}_{E_1}$, then
$$e+\alpha e=P^{E_*}_{E_0}P^{E_*}_{E_1}e+P^{E_0}_{E_*}P^{E_*}_{E_1}e=P^{E_*}_{E_1}e\in
E_1$$ for any $e\in E_0$; conversely,
$$e=P^{E_*}_{E_0}e+P^{E_0}_{E_*}e=P^{E_*}_{E_0}e+P^{E_0}_{E_*}P^{E_*}_{E_1}P^{E_*}_{E_0}e$$
for any $e\in E_1.$ So
$$E_1=\{e+\alpha e:\forall e\in E_0\}.$$
Conversely, assume $E_1=\{e_0+\alpha e_0:\forall e_0\in E_0\}$ for
any $\alpha\in B(E_0,E_*)$. We claim that $E_1$ is closed. Let
$e=e_0+\alpha e_0\rightarrow e_*$, then
$$P^{E_*}_{E_0}e=e_0\rightarrow P^{E_*}_{E_0}e_*\in E_0,{\rm and}\
 P^{E_0}_{E_*}e=\alpha
e_0\rightarrow\alpha(P^{E_*}_{E_0}e_*)$$
since for  all
$\alpha,P^{E_*}_{E_0}$ and $P^{E_0}_{E_*}$ are bounded. So
$$e_*=P^{E_*}_{E_0}e_*+P^{E_0}_{E_*}e_*=P^{E_*}_{E_0}e_*+\alpha(P^{E_*}_{E_0}e_*)\in
E_1.$$ This shows that $E_1$ is closed. Also, one observes $E_1\cap
E_*=\{0\}$, since $e_0+\alpha e_0\in E_*$ implies $e_0=0$.
Obviously,
$$e=P^{E_*}_{E_0}e+P^{E_0}_{E_*}e=(P^{E_*}_{E_0}e+\alpha(P^{E_*}_{E_0}e))+(P^{E_0}_{E_*}e-\alpha(P^{E_*}_{E_0}e)),\quad
\forall e\in E,$$ which shows $E_1\oplus E_*\supset E$. So
$E_1\oplus E_*=E$. The proof ends.\quad$\Box$

{\bf Theorem 1.11}\quad  {\it Suppose that $A\in B(E,F)$ is double
splitting and  $A\not=0$. Let  $A^+$ be a generalized inverse of
$A$. Define
$$M(A,A^+)(T)=(T-A)A^+A+C^{-1}_A(A^+,T)T,\quad\quad\forall T\in
V(A,A^+).$$ Then $M(A,A^+)$ is a smooth diffeomorphism from
$V(A,A^+)$ onto itself with the fixed point $A$.

For abbreviation, write $M(A,A^+)$ as $M$ in the sequel.}

{\bf Proof}\quad Note $C^{-1}_A(A^+,T)TA^+=AA^+$ as  indicated in
the proof of Theorem 1.1. Clearly,
\begin{eqnarray*}
(M(T)-A)A^+&=&(T-A)A^++C^{-1}_A(A^+,T)TA^+-AA^+\\
&=&(T-A)A^+\end{eqnarray*} so $C^{-1}_A(A^+,M(T))=C^{-1}_A(A^+,T).$
In addition, $C^{-1}_A(A^+,T)$ for $T\in V(A^+,A)$ is smooth, and
$C_A(A^+,A)=I_F$. Then one concludes $M(T)$ is a smooth map on
$V(A,A^+)$ into itself  with the fixed point $A$. Now we merely need
to show that $M$ has  an inverse map on $V(A,A^+)$. Fortunately, we
have
$$M^{-1}(m)=mA^+A+C_A(A^+,m)m(I_E-A^+A)\quad{\rm for\ all}\quad m\in
V(A,A^+).$$
In what follows, we are going to examine $(M\circ M^{-1})(m)=m$ for all $m\in V(A,A^+)$ and $(M^{-1}\circ M)(T)=T$  for all $T\in V(A,A^+)$.
Evidently,
$$(M^{-1}(m)-A)A^+=(m-A)A^++C_A(A^+, m)m(I_E-A^+A)A^+=(m-A)A^+.$$
Hereby, one observes   $C_A(A^+,M^{-1}(m))=C_A(A^+,m)$ and so
$$C^{-1}_A(A^+,m)=C^{-1}_A(A^+,M^{-1}(m)) \quad{\rm
for\ any}\quad  m\in V(A,A^+).$$
Thus
\begin{eqnarray*}
(M\circ M^{-1})(m)&=&(M^{-1}(m)-A)A^+A+C^{-1}_A(A^+,M^{-1}(m))M^{-1}(m)\\
&=&(m-A)A^+A+C^{-1}_A(A^+,m)(mA^+A+C_A(A^+,m)m(I_E-A^+A))\\
&=&(m-A)A^+A+m(I_E-A^+A)+A=m,\quad\forall m\in V(A,A^+);\\
(M^{-1}\circ M)(T)&=&M(T)A^+A+C_A(A^+,M(T))M(T)(I_E-A^+A)\\
&=&(T-A)A^+A+C^{-1}_A(A^+,T)TAA^+A+C_A(A^+,M(T))M(T)(I_E-A^+A)\\
&=&TA^+A+C_A(A^+,M(T))C^{-1}_A(A^+,T)T(I_E-A^+A)\\
&=&TA^+A+T(I_E-A^+A)=T,\quad\forall T\in V(A,A^+).\end{eqnarray*}
$\Box$

 \vskip 0.2cm\begin{center}{\bf 2\quad The co-final set and Frobenius
 Theorem}\end{center}\vskip 0.2cm

Let $E$ be a Banach space, and $\Lambda$ an open set in $E$. Assign
a subspace $M(x)$ in $E$ for every point $x$ in $\Lambda$,
especially, where the dimension of $M(x)$ may be infinite. In this section, we consider the
family $\cal{F}$ consisting of all $M(x)$ over $\Lambda.$ We
investigate the sufficient and necessary condition for $\cal{F}$
being $c^1$ integrable at  a pint $x$ in $\Lambda$. We stress  that
the concept of co-final set of $\cal{F}$ at $x_0\in\Lambda$ is
introduced.

{\bf Definition 2.1}\quad{\it Suppose $E=M(x_0)\oplus E_*$ for
$x_0\in\Lambda$. The set
$$J(x_0,E_*)=\{x\in\Lambda:M(x)\oplus E_*=E\},$$
is called a  co-final set of $\cal{F}$ at $x_0$.}

The following theorem tell us the co-final set and integral submanifold of $\cal{F}$ at $x_0$ has some connection.

{\bf Theorem 2.1}\quad{\it  If $\cal{F}$ is $c^1$-integrable at
$x_0\in\Lambda$, say that $S\subset E$ is an integral submanifold of
$\cal{F}$ at $x_0$, then there exist a closed subspace $E_*$ and a
neighborhood $U_0$ at $x_0$, such that}
$$M(x)\oplus E_*=E,\quad\quad \forall x\in S\cap U_0,$$
i.e., $$J(x_0,E_*)\supset S\cap U_0.$$

{\bf Proof}\quad Recall that the  submanifold $S$ in the Banach
space $E$ is said to be tangent to $M(x)$ at $x\in S$ provided
$$M(x)=\{\dot{c}(0):\forall c^1{\rm -curve}\subset S{\rm with}\ c(0)=x\}.\eqno(2.1)$$
By the definition of the submanifold, there exist a subspace $E_0$
splitting in $E$, say $E=E_0\oplus E_1$, a neighborhood $U_0$ at
$x_0$, and a $c^1$-diffeomorphism
$\varphi:U_0\rightarrow\varphi(U_0)$ such that $\varphi(S\cap U_0)$
is an open set in $E_0$. We claim
$$\varphi'(x)M(x)=E_0,\quad\quad\forall x\in S\cap U_0.\eqno(2.2)$$

Let $c(t)$ with $c(0)=x$ be an arbitrary $c^1$-curve contained in
$S\cap U_0$, then $\varphi'(x)\dot{c}(0)\in E_0$, and so
$\varphi'(x)M(x)\subset E_0$. Conversely, let
$r(t)=\varphi(x)+te\subset\varphi(S\cap U_0)$ for any $e\in E_0$,
and set $c(t)=\varphi^{-1}(r(t))\subset S\cap U_0$, then
$\varphi'(x)\dot{c}(0)=e$ and so $\varphi'(x)M(x)\supset E_0$. We
now conclude that (2.2) holds. For abbreviation in the sequel, write
$M_0=M(x_0)$. Let $E_x=\varphi'(x)^{-1}E_1$ and $E_*=E_{x_0}$ We
claim $M(x)\oplus E_*=E,\forall x\in S\cap U_0$. Evidently,
$\varphi'(x_0)(M_0\oplus E_*)=E_0\oplus E_1=E$ by (2.2), and so,
$M_0\oplus E_*=E$. Consider the projection
$$P_x=\varphi'(x)^{-1}\varphi'(x_0)P^{E_*}_{M_0}\varphi'(x_0)^{-1}\varphi'(x),\quad\quad\forall
x\in S\cap U_0.\eqno(2.3)$$ Obviously, $P^2_x=P_x$, i.e., $P_x$ is a
projection on $E$. Next go to show $R(P_x)=M(x)$ and $N(P_x)=E_x$.
Indeed,
\begin{eqnarray*}
e\in N(P_x)&\Leftrightarrow&\varphi'(x_0)^{-1}\varphi'(x)e\in
E_*\Leftrightarrow\varphi'(x)e\in E_1\\
&\Leftrightarrow&e\in\varphi'(x)^{-1}E_1=E_x;\\
R(P_x)&=&\varphi'(x)^{-1}\varphi'(x_0)M_0=\varphi'(x)^{-1}E_0=M(x).
\end{eqnarray*}
Obviously, $P_x$ is a generalized inverse of itself, and
$\lim\limits_{x\rightarrow x_0}P_x=P^{E_*}_{M_0}.$

Recall in the proof of the sufficiency of Theorem 1.5 the following conclusion is proved at first:
for some generalized inverse $T_0^+$ of $T_0$, $T_x$ near $x_0$ has a generalized inverse $T_x^+$
and $\lim\limits_{x\rightarrow x_0}T_x^+=T_0^+$, then there exists a neighborhood $V_0$ at $x_0$
such that $R(T_x)\cap N(T_0^+)=\{0\}$ for all $x\in V_0$. We mention that this conclusion
is very useful in the modern perturbation analysis of generalized inverses.

Now, take $S\cap U_0, P_x, P_{M_0}^{E_*}$ and $P_{M_0}^{E_*}$ in place of $X, T_x, T_0$ and $T_0^+$ in the conclusion, respectively.
Then we have $R(P_x)\cap N(P_{M_0}^{E_*})=\{0\}$ for all $x\in V_0$. Without loss of generality, one can assume
$V_0\subset\{x\in S\cap U_0: \|P_x-P_{M_0}^{E_*}\|<\|P_{M_0}^{E_*}\|^{-1} \}$ and $V_0=S\cap U_0$.
Thus, by the equivalence of the conditions (i) and (iii) in Theorem 1.1 we conclude
$$R(P_x)\oplus N(P_{M_0}^{E_*})=E, \, \text{i.e.},\,\, M(x)\oplus E_*=E, \,\,\, \forall x\in S\cap U_0.$$
This shows $J(x_0,
E_*)\supset S\cap U_0.$\quad$\Box$

Using $M_0$ and $E_*$ in place of $E_0$ and $E_1$ in Theorem 2.1
respectively, we have

{\bf Lemma 2.1}\quad{\it  With the same assumption and notations,
$E_*, E_0$ and $U_0$ as in Theorem 2.1, the following conclusion
holds: there exists a $c^1$ diffeomorphism $\varphi:U_0\rightarrow
\varphi(U_0)$ with $\varphi(x_0)=0$ such that $V_0=\varphi(S\cap
U_0)$ is an open set in $M_0$, and}
$$\varphi'(x)M(x)=M_0,\quad\quad\forall x\in S\cap U_0.\eqno(2.4)$$

{\bf Proof}\quad Consider the $c^1$-diffeomorphism $\varphi$ in the
proof of Theorem 2.1 and no loss of generality, we may assume
$\varphi(x_0)=0$, since otherwise we may take $\varphi-\varphi(x_0)$ in
place  of $\varphi$. Also, $\varphi'(x)M(x)=E_0$ for all $x\in S\cap
U_0$, specially, $\varphi'(x_0)M_0=E_0$. Again instead of $\varphi$
by $\varphi'(x_0)^{-1}\varphi$, still write it as $\varphi,$ the
lemma follows.\quad $\Box$

By Theorem 2.1, $J(x_0,E_*)\supset S\cap U_0$ when $\cal{F}$ is
integrable at $x_0$. Due to the co-final set $J(x_0,E_*)$, we get
the applicable equation of the integral submanifold under the
coordinate system $(M_0,0,E_*)$ as shown in the next lemma.

{\bf Lemma 2.2}\quad{\it  Suppose that $\cal{F}$ is $c^1$-integrable
at $x_0$. With the same notations, $S,U_0,V_0$ and $\varphi$ as in
Lemma 2.1, the following conclusion holds: there exists a
neighborhood $V_*$ at 0 in $V_0$, such that
$P^{E_*}_{M_0}\varphi^{-1}:V_*\rightarrow
P^{E_*}_{M_0}\varphi^{-1}(V_*)$ is a $c^1$ diffeomorphism in $M_0$, and
$$\varphi^{-1}(V_*)=(I_E+\psi)(V),\eqno(2.5)$$
where
$V=P^{E_*}_{M_0}\varphi^{-1}(V_*)$ and
$\psi=(P^{M_0}_{E_*}\varphi^{-1})\circ(P^{E_*}_{M_0}\varphi^{-1}|_{V_*})^{-1}.$}

{\bf Proof}\quad By Lemma 2.1, there are a neighborhood $U_0$ at
$x_0$ and $c^1$ diffeomorphism $\varphi$ from $U_0$ onto
$\varphi(U_0)$ with $\varphi(x_0)=0$, such that $V_0=\varphi(S\cap
U_0)$ is an open set in $M_0$. Consider the $c^1$ map
$\varphi_1=p^{E_*}_{M_0}\varphi^{-1}:V_0\rightarrow M_0$. In order
to seek $V_*$ in (2.5), we are going to show that $\varphi'_1(0)$ is
invertible in $B(M_0)$. According to (2.4) we have
$M_0=\varphi'(x_0)^{-1}M_0=(\varphi^{-1})'(0)M_0$, so that
$\varphi'_1(0)M_0=P^{E_*}_{M_0}(\varphi^{-1})'(0)M_0=M_0$; while
$$\varphi'_1(0)e=P^{E_*}_{M_0}(\varphi^{-1})'(0)e=(\varphi^{-1})'(0)e\quad{\rm
for\ any}\quad  e\in M_0,$$ and so $N(\varphi'_1(0))=\{0\}$ because
of $\varphi$ being a $c^1$ diffeomorphism. Hence $\varphi'_1(0)$ is
invertible in $B(M_0)$. Thus by the inverse map theorem, there
exists a neighborhood $V_*$ at 0 in $V_0$ such that
$\varphi_1:V_*\rightarrow\varphi_1(V_*)$ is a $c^1$ diffeomorphism.
Let $V=\varphi_1(V_*)$  and $\psi=\varphi_2\circ\varphi^{-1}_1$
where $\varphi_2=P^{M_0}_{E_*}\varphi^{-1}$. Then
$$\varphi^{-1}(V_*)=\varphi_1(V_*)+\varphi_2(V_*)=(I_E+\psi)(V).$$
The proof ends.\quad$\Box$

According to Theorem 1.10,  for $x\in J(x_0,E_*)$ $M(x)$ has the
coordinate expression
$$M(x)=\{e+\alpha (x)e:\forall e\in M_0\}$$
where $\alpha\in B(M_0,E_*).$

We now state the  Frobenius theorem in a Banach space.

{\bf Theorem 2.2}\ (Frobenius theorem)\quad{\it $\cal{F}$ is $c^1$
integrable at $x_0$ if and only if the following conditions hold:}

(i) $M_0$ {\it splits in} $E$, {\it say} $E=M_0\oplus E_*$;

(ii) {\it there exists a neighborhood} $V$ $at P^{E_*}_{M_0}x_0$
{\it in} $M_0$, {\it and a} $c^1$ {\it map} $\psi:V\rightarrow E_*$,
{\it such that} $x+\psi(x)\in J(x_0,E_*)$ {\it for all} $x\in V$,
{\it and} $\alpha$ {\it is continuous in} $V$;

(iii) $\psi$  satisfies
$$\begin{array}{rl}
&\psi'(x)=\alpha(x+\psi(x))\quad{\rm for\ all}\quad x\in V,\\
&\psi(P^{E_*}_{M_0}x_0)=P^{M_0}_{E_*}x_0,\end{array}\eqno(2.6)$$
{\it where} $\psi'(x)${\it  is Fr\'{e}chet derivative of} $\psi$ $at$
$x$.

{\bf Proof}\quad Assume that $\cal{F}$ is $c^1$-integrable at $x_0$.
Go to prove that the conditions (i), (ii) and (iii) in the theorem
hold. By Theorem 2.1, the condition (i) holds, and
$J(x_0,E_*)\supset S\cap U_0$. By Lemma 2.2, there is a neighborhood
$V_*$ at $0$ in $V_0=\varphi(S\cap U_0)$ such that
$$ J(x_0,E_*)\supset S\cap U_0\supset
\varphi^{-1}(V_*)=\{x+\psi(x):\forall x\in V\},$$
where
$$V=P^{E_*}_{M_0}\varphi^{-1}(V_*)\quad{\rm and}\quad
\psi=(P^{M_0}_{E_*}\varphi^{-1})\circ(P^{E_*}_{M_0}\varphi^{-1}|_{V_*})^{-1}.$$
Obviously, the equation (2.6) implies that $\alpha(x+\psi(x))$ is
continuous in $V$. Next we merely need to prove that the equation (2.6) holds. First we claim
$$M(x+\psi(x))=(I_E+\psi'(x))M_0\quad\quad\forall x\in
V.\eqno(2.7)$$
By (2.5),
$$d(t)=(I_E+\psi)(c(t))\quad{\rm for\ any}\quad c^1{\rm {\tiny-}curve}\ c(t)\subset V\ {\rm with}\ c(0)=x$$
is a  $c^1$-curve$\subset\varphi^{-1}(V_*)$ with $d(0)=x+\psi(x)$.
While, since $V_{*}\subset V_{0}$, $\varphi^{-1}(V_*)$ is the $c^1$ integral submanifold
of $\cal{F}$ at $x+\psi(x)$, and so, $M(x+\psi(x))$ is tangent to
$\varphi^{-1}(V_*)$ at $x+\psi(x)$. Then, one can conclude
$M(x+\psi(x))\supset(I_E+\psi')(x)M_0$ for any $x\in V.$ Conversely,
let
$$c(t)=P^{E_*}_{M_0}d(t)\quad{\rm for\ any}\quad c^1{\rm{\tiny-}curve}
\ d(t)\subset\varphi^{-1}(V_*)\ {\rm with}\ d(0)=x+\psi(x),$$ then
it follows that $c(t)$ is a $c^1$-curve$\subset V$ with $c(0)=x$
since
$$c(t)=P^{E_*}_{M_0}\varphi^{-1}(\varphi(d(t))),\quad\varphi(d(t))\subset
V_*,\quad{\rm and}\quad V=P^{E_*}_{M_0}\varphi^{-1}(V_*).$$ Note
$$\psi=(P^{M_0}_{E_*}\varphi^{-1})\circ(P^{E_*}_{M_0}\varphi^{-1}|_{V_*})^{-1}.$$
Evidently
\begin{eqnarray*}
(I_E+\psi)(c(t))&=&P^{E_*}_{M_0}d(t)+\psi(P^{E_*}_{M_0}d(t))\\
&=&P^{E_*}_{M_0}\varphi^{-1}(\varphi(d(t)))+\psi(P^{E_*}_{M_0}\varphi^{-1}(\varphi(d(t)))\\
&=&P^{E_*}_{M_0}\varphi^{-1}(\varphi(d(t))+P^{M_0}_{E_*}\varphi^{-1}(d(t)))\\
&=&P^{E_*}_{M_0}d(t)+P^{M_0}_{E_*}d(t)=d(t).
\end{eqnarray*}
So (2.7) holds.

In what follows, we go to verify (2.6). Due to the co-final set
$J(x_0,E_*)\supset\{x+\psi(x):\forall x\in V\}$ as pointed at the
beginning of the proof and Theorem 1.10, we have
$$M(x+\psi(x))=\{e+\alpha(x+\psi(x))e:\forall e\in M_0\},$$
where $\alpha(x+\psi(x))\in B(M_0,E_*).$ Thus, it follows from the
equality (2.7) that for any $e\in M_0$, there exists $e_0\in M_0$
such that
$$e+\alpha(x+\psi(x))e=e_0+\psi'(x)e_0.$$
Hence $e=e_0$, so that
$$\alpha(x+\psi(x))e=\psi'(x)e\quad{\rm  for \ all}\quad e\in M_0,$$
i.e.,
$$\alpha(x+\psi(x))=\psi'(x),\quad\quad\forall x\in V.$$
Obviously,
$$\psi(P^{E_*}_{M_0}x_0)=\psi(P^{E_*}_{M_0}\varphi^{-1}(0))=P^{M_0}_{E_*}(0)x_0.$$
Now, the necessity of the theorem is proved.

Assume that the conditions (i), (ii) and (iii) hold. Go to show that
$\cal{F}$ is $c^1$ integrable at $x_0$. Let $S=\{x+\psi(x):\forall
x\in V\}$,
$$V^*=\{x\in E:P^{E_*}_{M_0}x\in V\},\quad{\rm and}\quad
\Phi(x)=x+\psi(P^{E_*}_{M_0}x),\forall x\in V^*.$$ Obviously, $V^*$
is an open set in $E$, $\Phi(V)=S,$ and $\Phi$ is a $c^1$ map.
Moreover, we are going to prove that $\Phi:V^*\rightarrow\Phi(V^*)$
is a diffeomorphism. Evidently, if $\Phi(x_1)=\Phi(x_2)$, for
$x_1,x_2\in V^*$, i.e.,
$$P^{E_*}_{M_0}(x_1-x_2)+\psi(P^{E_*}_{M_0}x_1)-\psi(P^{E_*}_{M_0}x_2)+P^{M_0}_{E_*}(x_1-x_2)=0$$
then $P^{E_*}_{M_0}x_1=P^{E_*}_{M_0}x_2$, and so
$P^{M_0}_{E_*}x_1=P^{M_0}_{E_*}x_2$. This shows that
$\Phi:V^*\rightarrow \Phi(V^*)$ is one-to-one. Now, in order to show
that $\Phi:V^*\rightarrow\Phi(V^*)$ is a $c^1$ diffeomorphism, we
need only to show that $\Phi(V^*)$ is an open set in $E$. By the inverse
map theorem it is enough to examine that $\Phi'(x)$ for any $x\in
V^*$ is invertible in $B(E)$. According to the condition (iii)
$$\begin{array}{rllr}
\Phi'(x)&=&I_E+\psi'(P^{E_*}_{M_0}x)P^{E_*}_{M_0}\\
&=&P^{E_*}_{M_0}+\alpha(P^{E_*}_{M_0}x+\psi(P^{E_*}_{M_0}x))P^{E_*}_{M_0}+P^{M_0}_{E_*}\end{array}\eqno(2.8)$$
for any $x\in V^*$. If $\Phi'(x)e=0$, i.e.,
$P^{E_*}_{M_0}e+\psi'(P^{E_*}_{M_0}x)P^{E_*}_{M_0}e+P^{M_0}_{E_*}e=0,$
then $P^{E_*}_{M_0}e=0$ and so, $P^{M_0}_{E_*}e=0$. This says
$N(\Phi'(x))=\{0\}$ for any $x\in V^*$. Next go to  verify that
$\Phi'(x)$ is surjective.

For abbreviation, write $M(y)=M_*$ and
$y=P^{E_*}_{M_0}x+\psi(P^{E_*}_{M_0}x)$ for any $x\in V^*$.
Obviously, $y\in S$ and by the assumption (ii), $y\in J(x_0,E_*)$,
i.e., $M_*\oplus E_*=E$. Hence the following conclusion  holds: there exists $e_0\in M_0$ such that
$P^{E_*}_{M_*}e=e_0+\alpha(y)e_0$ for any
$e\in E$. Set $e_*=e_0+P^{M_*}_{E_*}e$,
then by (2.8),
$$\Phi'(x)e_*=e_0+\alpha(y)e_0+P^{M_*}_{E_*}e=P^{E_*}_{M_*}e+P^{M_*}_{E_*}e=e.$$

 This says that $\Phi'(x)$ for any $x\in V^*$ is surjective. Thus we
 have
 proved that $\Phi^{-1}$  is a
 $c^1$-diffeomorphism  from open set $\Phi(V^*)$ onto $V^*$ and $\Phi^{-1}(S)=V$ is an open set in $M_0$.
  i.e., $S$ is a $c^1$ submanifold of $E$. Finally go to show that
  $S$ is tangent to $M(x)$ at any point $x\in S$.

  Write
  \begin{center}$T(x+\psi(x))=\{\dot{c}(0):\forall c^1$-curve$c(t)\subset
  S\ {\rm with}\ c(0)=x+\psi(x)\}$\end{center}
  for any $x\in V$.

  Repeat the process of the proof of the equality (2.2) in Theorem
  2.1, one can conclude
  $$(\Phi^{-1})'(x+\psi(x))T(x+\psi(x))=M_0,\quad{\rm i.\
  e.,}\quad \Phi'(x)M_0=T(x+\psi(x))$$
  for any $x\in V$. By (2.8),
  $$\Phi'(x+\psi(x))e=e+\alpha(x+\psi(x))e,\quad\quad\forall e\in
  M_0.$$
  So
  $$T(x+\psi(x))=\Phi'(x)M_0=M(x+\psi(x)),\quad\quad\forall x\in V,$$
  which shows that $S$ is tangent to $M(x)$ at any $x\in S$.
  The proof ends.\quad$\Box$

  The co-final set $J(x_0,E_*)$ is essential to  Frobenius theorem in
  Banach space. When $J(x_0,E_*)$ is trivial, which means that $x_0$
  is an inner point of $J(x_0,E_*)$, the theorem reduces to the initial value problem  (2.6). The following
  example will illustrate this fact, although it is very simple.

  {\bf Example}\quad{\it Let $E=\mathbf{R}^2,
  \Lambda=\mathbf{R}^2\setminus(0,0)$ and
  $$M(x,y)=\{(X,Y)\in\mathbf{R}^2:Xx+Yy=0\},\quad\quad\forall
  (x,y)\in\Lambda.$$
  Consider the family of subspaces, $\mathcal{F}=\{M(x,y):\forall
  (x,y)\in \Lambda\}$. Applying Frobenius theorem in Banach space to
  determine the  integral curve of $\mathcal{F}$ at $(0,1)$.}

  Set $U_0=\{(x,y)\in\mathbf{R}^2:y>0\}$ and
  $E_*=\{(0,y)\in\mathbf{R}^2:\forall y\in R\}$.
  Obviously, $U_0\subset \Lambda$, and
  $$M(x,y)\oplus E_*=\mathbf{R}^2,\quad\quad\forall (x,y)\in U_0,$$
  since
  $$M(x,y)\cap E_*=(0,0),\quad\quad\forall(x,y)\in U_0.$$

  Hence $J((0,1),E_*)\supset U_0$. This shows that $J((0,1),E_*)$
  is trivial.

  Next go to determine $\alpha$ in the equation (2.6). Note
  $M_0=M(0,1)=\{(X,0):\forall X\in R\}$. Evidently,
  \begin{eqnarray*}
  M(x,y)&=&\{(X,-\frac{x}{y}X):\forall X\in R\}\\
  &=&\{(X,0)+(0-\frac{x}{y}X):\forall X\in
  R\}\end{eqnarray*}
  for all $(x,y)\in U_0$.
   By Theorem 1.10 we can conclude
  $$\alpha(x,y)(X,0)=(0,-\frac{x}{y}X),\quad\quad\forall (x,y)\in
  U_0.$$
  Moreover, we claim that $\alpha:U_0\rightarrow B(M_0,E_*)$ is
  continuous in $U_0$. Obviously
  \begin{eqnarray*}
  &&\left\|(\alpha(x+\Delta x,y+\Delta
  y)-\alpha(x,y))(X,0)\right\|\\
  &&\quad\quad=\left\|(0, (\frac{x}{y}-\frac{x+\Delta
  x}{y+\Delta y})X)\right\|=\left|\frac{x+\Delta x}{y+\Delta
  y}-\frac{x}{y}\right|\|(X,0)\|
  \end{eqnarray*}
  for any $(x,y)\in U_0$, where $\|,\|$ denotes the norm in
  $\mathbf{R}^2$.

  So
  $$\|\alpha(x+\Delta x,y+\Delta
  y)-\alpha(x,y)\|=\left|\frac{x+\Delta x}{y+\Delta
  y}-\frac{x}{y}\right|,$$
  where $\|,\|$ denotes the norm in $B(M_0,E_*).$ Hereby, one can
  conclude that $\alpha:U_0\rightarrow B(M_0,E_*)$ is continuous in
  $U_0$. Let $V=\{(x,0):|x|<1\}\subset M((0,1))$ and $\psi((x,0))=(0,y(x))\subset E_*$ for any
  $(x,0)\in V$, where $y(x):(-1, 1)\mapsto {\bf R}$ is a $c^1$ function. Then by the equation (2.6).
 \begin{eqnarray*}
&&\frac{dy}{dx}=-\frac{x}{y},\quad\quad\forall x\in (-1,1),\\
&&y(0)=1.\end{eqnarray*}
It is easy to see that the solution is
$y=\sqrt{1-x^2}$ for all $x\in(-1,1)$.
So, $S=\{(x, \sqrt{1-x^2}): x\in(-1, 1)\}$ is the smooth integral curve of $\mathcal F=\{M(x, y): (x, y)\in\Lambda\}$ at $(0,1)$.
(See [Ma4] also.)

 When the co-final set $J(x_0,E_*)$ of $\mathcal{F}$ at $x_0$ is
 non trivial, it is a key point to Frobenius theorem to seek $J(x_0,E_*)$.
 In the next section, we will consider such a family of subspaces,
 which appears in the investigation of geometrical method for some
 partial differential equations (see [Caf]).

\vskip 0.2cm\begin{center}{\bf 3\quad A Family of Subspaces with
Non-trivial Co-Final Set and Its Smooth Integral
Submanifolds}\end{center}\vskip 0.2cm

Let  $\Lambda =B(E,F)\setminus\{0\}$ and $M(X)=\{T\in
B(E,F):TN(X)\subset R(X)\}$ for $X\in\Lambda$. V. Cafagna introduced
the  geometrical method for some partial differential equations and
presented the family of subspaces $\mathcal{F}=\{M(X)\}_{X\in\Lambda}$
in [Caf]. Now we take the following example to  illustrate that when
A is neither left nor right invertible in $B(R^2), $ the co-final
set of $\cal{F}$ at $A$ is non-trivial.

{\bf Example}\quad{\it  Let $\Lambda =B(R^2)\setminus\{0\}$,  the
matrix $A=\{a_{i,j}\}^2_{i,j=1}$ and the subspace $\mathbb{E}_*$ in
$B(R^2)$ be as follows
$$a_{i,j}=0\quad{\rm except}\quad a_{1,1}=1,\quad{\rm and}\quad
\{T=\{t_{i,j}\}\in B(R^2):t_{i,j}=0\quad{\rm except}\quad
t_{2,2}\},$$
 respectively. To verify
that $J(A,\mathbb{E}_*)$ is non-trivial. Evidently,
$N(A)=\{(0,y):\forall y\in R\}$ and $R(A)=\{(x,0):\forall x\in R\}$.
So
$$M(A)=\{T=\{t_{i,j}\}\in B(R^2):t_{2,2}=0\}$$
and
$$M(A)\oplus \mathbb{E}_*=B(\mathbf{R}^2).$$
However, there is an invertible matrix $A_\varepsilon=\{a_{i,j}\}$
with $a_{1,1}=1, a_{1,2}=a_{2,1}=0$, and $a_{2,2}=\varepsilon$,
 such that $M(A_\varepsilon)=B(R^2)$ for any $\varepsilon\not=0$.
 This says that $J(A,\mathbf{E}_*)$ is non-trivial since
 dim$M(A)=3, $dim$M(A_\varepsilon)=4$, and $\lim\limits_{\varepsilon\rightarrow
 0}A_\varepsilon=A$.}

  In this case, for seeking integral submanifold of $\mathcal{F}$
  at $A$, we try the co-final set $J(A,\mathbb{E}_*)$ according to
  Theorem 2.1.

  {\bf Lemma 3.1}\quad{\it Suppose that $X\in\Lambda$ is double
  splitting, say that  $X^+$ is a generalized inverse of $X$. Then
  $$M(X)=\left\{P^{N(X^+)}_{R(X)}T+R^{R(X)}_{N(X^+)}TP^{N(X)}_{R(X^+)}:\forall
  T\in B(E,F)\right\},\eqno(3.1)$$
  and one of its complements is}
  $$\mathbb{E}_X=\left\{P^{R(X)}_{N(X^+)}TP^{R(X^+)}_{N(X)}:\forall
  T\in B(E,F)\right\}.\eqno(3.2)$$

  {\bf Proof}\quad Obviously
  $$T=P^{N(X^+)}_{R(X)}T+P^{R(X)}_{N(X^+)}TP^{N(X)}_{R(X^+)}+P^{R(X)}_{N(X^+)}TP^{R(X^+)}_{N(X)}.\eqno(3.3)$$
  So the equality (3.1) follows from the definition of $M(X)$ and
  (3.3).\quad$\Box$

  As a corollary of the lemma we have

  {\bf Corollary 3.1}\quad {\it With the assumption and notations as in
  Lemma 3.1 we have}
  $$\mathbb{E}_X=\left\{T\in B(E,F):R(T)\subset N(X^+)\quad{\rm
  and}\quad N(T)\subset R(X^+)\right\}.\eqno(3.4)$$

  This is immediate from (3.1) and (3.2).

  Let $\mathbb{E}_*=\{T\in B(E,F):R(T)\subset N(A^+)$ and $N(T)\subset R(A^+)\}$.
  By Lemma 3.1 and Corollary 3.1, $M(A)\oplus \mathbb{E}_*=B(E,F).$
  Next, we go to determine the integral submanifold of $\mathcal{F}$ at $A$.
  Let
  $$S = \{X\in V(A, A^{+}) : R(X)\cap N(A^{+}) = \{0\} \}.$$

  Since the conditions (i), (ii) and (iii) in Theorem 1.1  are equivalent each other,
  $S$ is the set of all double splitting operators $X \in V(A,A^{+})$
  with $R(X^{+}) = R(A^{+})$ and $N(X^{+}) = N(A^{+})$. Then by Corollary 3.1,
  $\mathbb{E}_{X} = \mathbb{E}_*$ for any $X\in S$. Therefore
  $$J(A,\mathbb{E}_*)\supset S=\{T\in V(A,A^+):R(T^+)=R(A^+)\ {\rm and}\
   N(T^+)=N(A^+)\}.$$

   Let $X\in\Lambda$ be double splitting. Consider the smooth
   diffeomorphism from $V(X,X^+)$ onto itself with the fixed point $X$,
   $$M(X,X^+)(T)=(T-X)X^+X+C^{-1}_X(X^+,T)T$$
   as defined in Theorem 1.11.

   By direct computing,
   $$\begin{array}{rl}
   &M(X,X^+)'(T)\Delta T\\
   &\quad=\Delta TX^+X+C^{-1}_X(X^+,T)\Delta
   T-C^{-1}_X(X^+,T)\Delta TX^+C^{-1}_X(X^+,T)T\end{array}
   \eqno(3.5)$$
   for any $T\in V(X,X^+)$, where $M(X,X^+)'(T)$ denotes the Fr$\acute{\rm e}$chet
   derivative  of $M(X,X^+)(T)$ at $T$. Moreover, we have

   {\bf Theorem 3.1}\quad{\it Suppose that $A\in\Lambda$ is double
   splitting. Then $S$ is a smooth submanifold in $B(E,F)$ and tangent
   to $M(X)$ at any $X\in S$.}

   {\bf Proof}\quad  First we show that $S$ is a smooth submanifold in $B(E,F)$. By Theorem 1.11,
    $M(A,A^+)(T)$ is a smooth
   diffeomorphism from $V(A,A^+)$ onto itself with the fixed point $A$,
   and hence we merely  claim
   $$M(A,A^+)(S)=M(A)\cap V(A,A^+).\eqno(3.6)$$
   Evidently,
   $$M(A,A^+)(T)N(A)=C^{-1}_A(A^+,T)TN(A)\subset R(A)$$
   because of the equivalence of the conditions (v) and (i) in Theorem
   1.1. This shows $M(A,A^+)(S)\subset M(A)\cap V(A,A^+).$
   Conversely, go to verify
   $$m=M(A,A^+)^{-1}(T)\in S\quad\forall T\in M(A)\cap V(A,A^+).$$
Note the equality shown in the proof of Theorem 1.1,
$C^{-1}_A(A^+,T)=C^{-1}_A(A^+,m).$  Let $m=M(A,A^+)^{-1}(T)$, for
any $T\in M(A)\cap V(A,A^+).$   Then
   $$C^{-1}_A(A^+,T)=C^{-1}_A(A^+,M(A,A^+)(m))=C^{-1}_A(A^+,m).$$
   Hereby,
   \begin{eqnarray*}
   C^{-1}_A(A^+,m)m N(A)&=& C^{-1}_A(A^+,T)M(A,A^+)^{-1}(T)N(A)\\
   &=&TN(A)\subset R(A),
   \end{eqnarray*}
   so that $m\in S$ because of the equivalence of (v) and (i) in Theorem
   1.1. This says   that  $S$ is a smooth submanifold in $B(E,F)$.
   Finally go  to  show that  $S$ is tangent to $M(X)$ at any
   $X\in S$. Clearly, $C_X(X^+,X)=I_F$ and so, $C^{-1}_X(X^+,X)=I_F$.
   Hereby
   $$M(X,X^+)'(X)\Delta T=\Delta T,\quad\quad\forall\Delta T\in
   B(E,F)\eqno(3.7)$$
   according to (3.5). First, we claim that $S$ is tangent to $M(A)$ at $A$.
   From (3.6) and (3.7) it follows that
   $M(A,A^+)'(A)\dot{c}(0)=\dot{c}(0)\in M(A)$ for any $c^1$-curve
   $c(t)\subset S$ with $c(0)=A$. This shows that the set of all tangent
   vectors of $S$ is contained in $M(A)$. Conversely, let
   $c(t)=M(A,A^+)^{-1}(d(t))$ for any $c^1$-curve $d(t)\subset
   M(A)\cap V(A,A^+)$ with $d(0)=A$. By (3.6), $c(t)\subset S$ and
   $c(0)=A$. While $d(t)=M(A,A^+)(c(t))$ and so, $\dot
   {d}(0)=M(A,A^+)'(A)\dot{c}(0)=\dot{c}(0).$ This shows that $M(A)$
   is contained in the set of all tangent vectors of $S$. Now it is
   proved that $S$ is tangent to $M(A)$ at $A$.
    Next turn to  the proof of the conclusion for any
    $X\in S$.
     As shown in the above, $X$ in $S$ has the the generalized inverse
     $X^{\oplus}$ with
    $R(X^\oplus)=R(A^+)$ and $N(X^\oplus)=N(A^+)$. So we can consider
    $S_1=\{T\in V(X, X^\oplus):R(T)\cap N(X^\oplus)=\{0\}\}$ for $X$ in $S$.
    For simplicity, we still write  $X^{\oplus}$ as $X^+$. Clearly,
    $$S\cap V(X, X^+)=S_1\cap V(A, A^+)=S_1\cap S \,\,\text{and}\,\, X\in S\cap S_1.$$
   In the same way  as the proof
   of (3.6) we have
   $$M(X,X^+)(S\cap S_1)=M(X)\cap V(A,A^+)\cap V(X,X^+).\eqno(3.8)$$
   Repeat the process of the proof of $S$ being tangent to $M(A)$ at
   $A$. Let $c(t)$ be any  $c^1$ curve$\subset S\cap S_1$ with
   $c(0)=X$. Then it follows from (3.8) and (3.7) that
   $$M(X,X^+)'(X)\dot{c}(0)=\dot{c}(0)\quad{\rm so\ that}\quad
   \dot{c}(0)\in M(X).$$
   Conversely, let $d(t)$ be any $c^1$-curve$\subset M(X)\cap
   V(A,A^+)\cap V(X,X^+)$ with $d(0)=X$ and
   $c(t)=M(X,X^+)^{-1}(d(t)).$ Then $d(t)=M(X,X^+)(c(t))$. By (3.7),
   $$\dot{d}(0)=M(X,X^+)'(X)\dot{c}(0)=\dot{c}(0).$$

   Combining the two results above one concludes that $S$ is tangent
   to $M(X)$ at $X$.  \quad $\Box$

{\bf Theorem 3.2}\quad{\it  Each of $F_k,\Phi_{m,n},\Phi_{m,\infty}$
  and $\Phi_{\infty,n}$ is a smooth submanifold in $B(E,F)$ and
  tangent to $M(X)$ at any $X$ in it.}

  {\bf Proof}\quad It is well known that $X\in
  F_k\cup\Phi_{m,n}\cup\Phi_{m,\infty}\cup\Phi_{\infty,n}$ is double
  splitting. By Theorem 3.1, $S=\{T\in V(X,X^+):R(T)\cap N(X^+)=\{0\}\}$
  is a
  smooth submanifold in $B(E,F)$  tangent to $M(T)$ at any $T\in
   S$, and $M(X,X^+)(S)=M(X)\cap V(X,X^+). $ Write $\Phi$ as any one
   of $F_k,\Phi_{m,n},\Phi_{m,\infty}$ and $\Phi_{\infty,n}$. By
   Theorem 1.2 and (3.6),
   $$S=\{T\in V(X,X^+):R(T)\cap N(X^+)=\{0\}\}\}=\Phi\cap V(X,X^+)$$
   and
   $$M(X,X^+)(S)=M(X,X^+)(\Phi\cap V(X,X^+))=M(X)\cap
   V(X,X^+).\eqno(3.9)$$
   Obviously, $(M(X,X^+),V(X,X^+),B(E,F))$ is a smooth admissible
   chart of $B(E,F)$ at $X\in\Phi$. While $M(X,X^+)(\Phi\cap
   V(X,X^+))=M(X)\cap V(X,X^+)$ is an open set in $M(X)$ because of
   (3.9) and $M(X)\oplus\mathbb{E}_X=B(E,F)$ by Lemma 3.1. So we
   conclude that $\Phi$ is a smooth submanifold in $B(E,F)$. Finally,
   by Theorem 3.1,$\Phi$ is tangent to $M(X)$ at any $X\in\Phi$. The
   theorem is proved.\quad$\Box$

   Theorem 3.2 expands the result for $\Phi_{1,1}$ to wider
   classes of operators. It seems to be nice to
   further developing of the method by V. Cafagna in [Caf].(see
   [An]).

\vskip 0.1cm
\begin{center}{\bf References}
\end{center}
\vskip -0.1cm
\medskip
{\footnotesize
\def\REF#1{\par\hangindent\parindent\indent\llap{#1\enspace}\ignorespaces}

\REF{[Abr]}\ R. Abraham, J. E. Marsden, and T. Ratin, Manifolds,
tensor analysis and applications, 2nd ed., Applied Mathematical
Sciences 75, Springer, New York, 1988.

\REF{[An]}\ V. I. Arnol'd, Geometrical methods in the theory of
ordinary differential equations, 2nd ed., Grundlehren der
Mathematischen Wissenschaften 250, Springer, New York, 1988.

\REF{[Beg]}\ M. S. Berger, Non-linearity and functional analysis, A
cademic Press, New York, 1977.

\REF{[Caf]}\ V. Cafagra, Global invertibility and finite
solvability", pp. 1-30 in Nonlinear functional analysis (Nework, NJ,
1987), edited by P. S. Milojevic, Lecture Notes in Pure and Appl.
Math. 121, Dekker, New York, 1990.

\REF{[H-M]}\ Qianglian Huang, Jipu Ma, Perturbation analysis of
generalized inverses of linear operators in Banach spaces, Linear
Algebra and its Appl., 389(2004), 359-364.

\REF{[Ka]}\ Kato, T., Perturbation Theory for Linear Operators, New
York: Springer-Verlag, 1982.

\REF{[Ma1]}\ Jipu Ma, (1,2) inverses of operators between Banach
spaces and local conjugacy theorem, Chinese Ann. Math. Ser. B.
20:1(1999), 57-62.

\REF{[Ma2]}\ Jipu Ma, A generalized preimage theorem in global
analysis, Sci, China Ser. A 44:33(2001), 299-303.

\REF{[Ma3]}\ Jipu Ma, A generalized transversility in global
analysis, Pacif. J. Math.,  236:2(2008),  357-371.

\REF{[Ma4]}\ Jipu Ma, A geometry characteristic of Banach spaces with $C^1$-norm, Front. Math. China,  2014, 9(5): 1089-1103.

 \REF{[N-C]}\ M. Z. Nashed, X. Chen, Convergence of Newton-like methods for singular equations using outer
 inverses, Numer
 Math., 66:(1993), 235-257.

 \REF{[P]}\ R. Penrose, A generalized inverse for Matrices, Proc. Cambridge Philos. Soc. 1955, 51:406-413.

 \REF{[Zei]}\ E. Zeidler,  Nonlinear functional analysis and its applications, IV: Applications to mathematical physics,
  Springer, New York, 1988.

  1. Department of Mathematics, Nanjing University, Nanjing, 210093,
  P. R. China

  2. Tseng, Yuanrong Functional Research Center, Harbin Normal
  University, Harbin, 150080, P. R. China

 E-mail address: jipuma@126.com
}
\end{document}